\documentclass[10pt]{amsart}

\usepackage{amsmath}
\usepackage{amsthm}
\usepackage{amsfonts}
\usepackage{amssymb,latexsym}
\usepackage[all]{xy}
\usepackage{graphicx}
\usepackage[mathscr]{eucal}
\usepackage{verbatim}
\usepackage{hyperref}

\theoremstyle{plain}
    \newtheorem{thm}{Theorem}[section]
    \newtheorem{prop}[thm]{Proposition}

    \newtheorem{cor}[thm]{Corollary}

\theoremstyle{definition}

\theoremstyle{remark}
    \newtheorem{rem}[thm]{Remark}

\numberwithin{equation}{section}

\begin{document}

\title{What is special about the Divisors of 24?}
\date{\today}

\author{Sunil K. Chebolu}
\address{Department of Mathematics \\
Illinois State University \\
Normal, IL 61790, USA} \email{schebol@ilstu.edu}

\maketitle

\emph{It is a miracle that the human mind can string a thousand arguments together without getting itself into contradictions.}
\begin{flushright}
-- Eugene Wigner
\end{flushright}

\thispagestyle{empty}
\section{Introduction}
The divisors of $24$ are $1, 2, 3, 4, 6, 8, 12$, and $24$. To pique the reader's interest  let us  pose the following riddle:
what is an \emph{interesting} number theoretic characterization of the divisors of 24 amongst all positive integers?   I will present one in terms of  modular multiplication tables.
This idea evolved  from a question raised by Elliott Mahler in my elementary number theory class.
Shortly after introducing the new world of $\mathbb{Z}_{n}$, I asked my students to write down the  multiplication tables for $\mathbb{Z}_{2}, \mathbb{Z}_{3}$, and $\mathbb{Z}_{4}$.  I then showed them the multiplication table of $\mathbb{Z}_{8}$  with the intention of drawing their attention to some differences between tables for prime and composite moduli.
\[\mathbb{Z}_8:  \ \
\begin {array}{c|cccccccc} * &0&1&2&3&4&5&6&7\\ \hline 0&0&0&0&0&0&0&0&0\\1&0&1&2&3&4&5&6&7\\ 2&0&2&4&6&0&2&4&6\\3 &0&3&6&1&4&7&2&5\\ 4&0&4&0&4&0&4&0&4\\ 5&0&5&2&7&4&1&6&3 \\ 6&0&6&4&2&0&6&4&2\\ 7&0&7&6&5&4&3&2&1
\end {array}
\]
Upon seeing these tables Elliott asked,
``\emph{I see that 1's in these multiplication tables appear only on the diagonal \footnote{Throughout this paper the term ``diagonal" refers to the main diagonal.}. Is that always true?}''
 Of course, looking further one knows that this is not always true. For instance, in table for  $\mathbb{Z}_{5}$,  $1$ occurs at an off-diagonal position $(2, 3)$   corresponding to the multiplication
\[ (2) (3) = 1   \ \text{ in } \mathbb{Z}_{5}. \]

Having seen some examples with 1's only on the diagonal and some with 1's also off the diagonal,  the following question begs to be answered.
\noindent
\emph{For what values of $n$ do 1's occur only on the diagonal in the multiplication table of $\mathbb{Z}_{n}$, never off the diagonal?}

I will investigate this question using various tools from number theory and will tie it up with some interesting topics which seem a priori unrelated to this question. Specifically, the tools used are:  the Chinese remainder theorem,  Dirichlet's theorem on primes in an arithmetic progression, the structure theory of units in $\mathbb{Z}_n$,  the Bertrand-Chebyshev theorem and its extension by  Erd\"{o}s and Ramanujan. There is no doubt that some of these tools are rather heavy-duty for the relatively simple question under investigation. However, my goal in this paper is to expose  students to as many interesting topics in number theory as possible via this question which came up naturally in the classroom, and to show them the interconnections between these various topics.

The question under investigation is answered by the following theorem.
\begin{thm} \label{main}
The multiplication table for $\mathbb{Z}_{n}$ contains  1's  only on the diagonal if and only if $n$ is a divisor of $24$.
\end{thm}
Note that the trivial  divisor $1$ of $24$  corresponds to $\mathbb{Z}_{1}$, which consists of only one element ($0$). Therefore the requirement of having all ones on the diagonal is vacuously satisfied in this case.

I will give five different arguments.  The first three provide a complete proof of Theorem \ref{main}, while the last two proofs show that the integers $n$ with the diagonal property satisfy $n \le 24$.  The finitely many values of $n$ up to $24$ can be dealt with separately to prove the main theorem. I will begin with a convenient proposition in the next section.

\vskip 3mm

\noindent
\textbf{Acknowledgements:} I would like to thank my student Elliott Mahler for raising the aforementioned question which led to this paper and also the anonymous referees for their comments and suggestions which helped improve the exposition. Finally, I want to thank my wife Surekha Methuku for her help with some MAPLE programs in connection to this paper.

\section{The diagonal condition}
Let me begin by   examining  the condition ``1's in the multiplication table for $\mathbb{Z}_{n}$ occur only along the diagonal'' more closely.
For convenience, I will refer to this as the diagonal condition/property  for $n$.  Let us  fix representatives for the elements in $\mathbb{Z}_n$:
\[ \mathbb{Z}_n = \{ 0,1,2, \cdots, n-1\}.\]
Suppose there is a 1 at position $(a, b)$ in the multiplication  table for $\mathbb{Z}_{n}$. This means $ab = 1$ in $\mathbb{Z}_{n}$. ($a$, and hence also $b$, is then said to be invertible in $\mathbb{Z}_n$.) If  the diagonal condition holds for $n$, then $(a, b)$ has to be a diagonal position. This would mean that  $a = b$, and therefore $a^{2}= 1$ in $\mathbb{Z}_{n}$, or equivalently $n$ divides $a^2-1$. It is an easy exercise to show that  $a$ is invertible in $\mathbb{Z}_{n}$  if and only if  $\gcd(a, n) = 1$.   We then have the following proposition.

\begin{prop} \label{prop}
Let $n$ be a positive integer. Then the following statements are equivalent.
\begin{enumerate}
\item $1$'s in the multiplication table for $\mathbb{Z}_n$ occur only on the diagonal.
\item If $a$ is an invertible element in $\mathbb{Z}_n$, then $a^2 = 1$ in $\mathbb{Z}_n$.
\item If $a$ is a positive integer that is relatively prime to $n$, then $n$ divides $a^2 - 1$.
\item If $p$ is a prime number that does not divide $n$, then $n$ divides $p^2 - 1$.
\end{enumerate}
\end{prop}

\begin{proof}
In light of the above discussion, the equivalence of the first three statements is clear. Moreover, (4) is a special case of (3). So it is enough to show that
(4) implies (3). To this end, let  $a$ be a positive integer that is relatively prime to $n$. If $a$ is 1, then the conclusion is obvious. If $a >1$, consider the
prime factorization $a = p_1^{a_1} p_2^{a_2} \cdots p_k^{a_k}$ of $a$. Since $a$ is relatively prime to $n$, none of these primes divide $n$. So by (4),
we have $p_i^2 \equiv 1 \mod{n}$ for all $i$. Then we have
\[ a^2 = (p_1^{a_1} p_2^{a_2} \cdots p_k^{a_k})^2 =  (p_1^2)^{a_1} (p_2^2)^{a_2} \cdots (p_k^2)^{a_k} \equiv (1)(1)\cdots(1) \equiv 1 \mod{n},\]
as desired.
\end{proof}

I will use the equivalent statements of Proposition \ref{prop} interchangeably when referring to integers which have the diagonal property.

\section{The Chinese remainder theorem}
In this section I will use the Chinese remainder theorem to give what is probably the shortest  proof of Theorem \ref{main}.
The Chinese remainder theorem in its classical form talks about simultaneous solutions to a system of linear congruences. It can be restated succinctly as an isomorphism of rings \cite[Page 265]{DummitFoote}
\[ \mathbb{Z}_{ab} \cong \mathbb{Z}_{a} \oplus \mathbb{Z}_b, \]
whenever $a$ and $b$ are positive integers that are relatively prime.
(Multiplication in $\mathbb{Z}_{a} \oplus \mathbb{Z}_b$ is done component-wise.) Isomorphism of rings means there is a 1-1 correspondence between  $\mathbb{Z}_{ab}$ and $ \mathbb{Z}_{a} \oplus \mathbb{Z}_b$ such that under this correspondence, the addition and multiplication tables in these two rings are the same. As an example, we ask the reader to verify that $\mathbb{Z}_{2} \oplus \mathbb{Z}_3 \cong \mathbb{Z}_6$ under the correspondence $(0, 0)\leftrightarrow 0, (0, 1)\leftrightarrow 4, (0, 2)\leftrightarrow 2, (1, 0)\leftrightarrow 3, (1, 1)\leftrightarrow 1,  (1, 2)\leftrightarrow 5$.

 To start, we first consider the case when $n$ is odd, so that $\gcd(2, n) = 1$. Then in order for $n$ to have the diagonal property,  $n$ has to  divide $2^{2} - 1 =  3$. This means $n$ has to be either $1$ or $3$, both of which have the diagonal property.
Next consider the case where $n$ is a power of 2, say $n = 2^{t}$ for some $t$, so that  $\gcd(3, n) = 1$. As before, for $n$ to have the diagonal property, $n$ has to divide $3^{2} - 1 = 8$. It is easily seen that all the divisors of $8$ have the diagonal property.
Now any positive integer $n$ can be uniquely written as
\[n = 2^t k,\]
where $k$ is odd and $t$ is a non-negative integer. Then by the Chinese remainder theorem  we have the isomorphism
\[ \mathbb{Z}_n \cong \mathbb{Z}_{2^t} \oplus \mathbb{Z}_k.\]
From this isomorphism it is easy to see that $n$ has the diagonal property if and only if both $2^t$ and $k$ have the diagonal property.
Combining these pieces, it follows that the only integers with the diagonal property are the divisors of $(8)(3) = 24$.

\section{Dirichlet's theorem on primes in an arithmetic progression}
Dirichlet proved the following theorem in 1837  which is a far-reaching extension of Euclid's theorem on the infinitude of primes. It states that given any two integers $s$ and $t$ that are  relatively prime, the arithmetic progression $\{ s x + t \, | \, x \text{ is a non-negative integer} \}$ contains infinitely many prime numbers; see \cite[Page 401]{NivenZuckerman}.  This result is one of the most beautiful results in all of number theory.

Let $n$ be an integer which has the diagonal property. So $n$ has the property that, for any prime $p$ which does not divide $n$,  $n\, |\, p^{2} - 1$. If $n \, | \, p^{2 } - 1$,  then for every prime divisor $q$ of $n$,  $q$ divides either $p -1$ or $p+1$.   In other words,  every prime $p$  that does not divide  $n$ has to be of the form $qx + 1$ or $qx -1$ for each prime divisor $q$ of $n$.  This is clearly a very strong condition on $n$.  If there is a prime divisor $q_{0}$ of $n$ which is bigger than $3$, then there will be an arithmetic progression
$\{ q_{0}x + r \, | \, x \ge 0\}$, where $r \ne 0, 1,$ or $q_0-1$ and $2 \le r \le q_{0}-2$. Note that $q_0$ and $r$ are then relatively prime, and
therefore Dirichlet's theorem  tells us that this arithmetic progression contains infinitely many primes.
 In particular, it contains a prime $p_{0}$ that does not divide $n$.  This choice of $p_{0}$ does not meet the requirement that  it is either of the form $q_{0}x + 1$ or $q_{0}x-1$. The upshot is that there is no prime divisor of $n$ which is bigger than  $3$, which means $n$ is of the form $2^{u}3^{v} $. The smallest prime number that is relatively prime to every number of the form $2^{u}3^{v}$ is $5$.  Our proposition then tells that $n$ has to divide $5^{2} - 1 = 24$, as desired.

\begin{rem} One can avoid the full strength of Dirichlet's theorem as used here. It is enough to assume the special case that the arithmetic
progression $5n + 2$ (or $5n+3$) contains infinitely many primes. This will allow us to show (exactly as above) that $5$ cannot divide $n$. Therefore it follows that  $n$ has to divide $5^2 -1 = 24$.  The above proof is, however,  more natural. It explains naturally why only primes $2$ and $3$ can occur in the factorization of $n$.
\end{rem}

\section{The structure theory of units in $\mathbb{Z}_n$}
The set of invertible elements (a.k.a. units) in $\mathbb{Z}_n$ is denoted by $U_n$. This set forms an abelian group under multiplication.  The structure of the group  $U_n$ has been completely determined. To explain, let $n = p_{1}^{c_{1}}p_{2}^{c_{2}} \cdots p_{k}^{c_{k}}$ be the prime decomposition of $n\, (>1)$.  The Chinese remainder theorem gives  an isomorphism of groups
\[ U_{n } \cong U_{p_{1}^{c_{1}}}  \oplus U_{p_{2}^{c_{2}}}\ \cdots  \oplus U_{p_{k}^{c_{k}}}. \]
(Again, isomorphism of groups means there is a 1-1 correspondence between the two groups such that the group multiplication tables under this correspondence are the same.)
It is therefore enough to explain the structure of $U_{p^{c}}$. This is given by \cite[Page 124]{NivenZuckerman}:
\[
U_{p^{c}} \cong
\begin{cases}
C_{1}   &  \text{ if }  \ \  p^{c} = 2^{1} \\
C_{2}   &   \text{ if } \ \ p^{c} = 2^{2} \\
C_{2} \oplus C_{2^{c-2}} & \text{ if } \ \ p^{c} = 2^{c} \text{ and } c \ge 3\\
C_{\phi(p^{c})} & \text{ if } \ \ p \text{ is odd},
\end{cases}
\]
where $C_{k}$ is the cyclic group of order $k$, and $\phi(x)$ is Euler's  function which denotes the number of positive integers less than $x$ that are relatively prime to $x$.

Returning to our problem, recall that  from  the above proposition, $n$ has the diagonal property if and only if  $a^{2} = 1$ for all $a$ in $U_n$.  Therefore our job is to simply identify those groups from the above list which have the property that every element  in them has order at most 2. $C_{1}$ and $C_{2}$ obviously  have this property. $C_{2} \oplus C_{2^{c-2}}$ will have this property if and only if $c - 2 \le 1$, or $c \le 3$. Finally, $C_{\phi(p^{c})}$ will have this property for $p$ is odd if and only if  $\phi(p^{c}) = p^{c-1}(p-1) \le 2$. It is easy to
see that this last inequality holds only  when $p^{c} = 3$.  From these calculations, we note that  an integer $n$ with the diagonal property cannot have a prime divisor bigger than $3$. Moreover, the maximum power of $3$ in $n$ has to be $1$, and that of $2$ has to be $3$.  The collection of these  integers is given by
\[ n = 2^{u}3^{v}, \ \ \text{where} \ \  0 \le u \le 3, \ \ 0 \le v \le 1,\]
which are exactly the divisors of $24$.

\begin{rem}  Note that the abelian group $U_{n}$ has a natural $\mathbb{F}_{2}$-vector space structure  precisely when $a^{2} = 1$ for all
$a$ in $U_{n}$.  Therefore we can say that $n$ has the diagonal property if and only if $U_{n}$ is naturally a vector space over $\mathbb{F}_{2}$.
\end{rem}

In the next two sections I will use some results in number theory to show that if a positive integer $n$ has the diagonal property then $n \le 24$.  The finitely many values of $n$ up to $24$ can then be dealt with separately to prove the main theorem.

\section{The Bertrand-Chebyshev theorem}
In the year 1845, Bertrand postulated that if $n \ge 2$ then there is always a prime number $p$ such that $n < p < 2n$. Although he did not give a proof, he
verified it for all values of $n$ up to three million. A few years later (1852) Chebyshev gave an  analytical proof of this result. Elementary proofs, however, had to wait until the next century.   In 1919 Ramanujan \cite{Ram1919} gave the first elementary proof using some properties of the gamma function and the Stirling's formula. His proof could be easily presented without ever mentioning the gamma function. In his first paper in 1932, Erd\"{o}s  \cite{Erdos}  gave another elementary proof of this theorem using some properties of the binomial coefficients.  See \cite{Book} for a nice presentation of Erd\"{o}s's proof. Hardy and Wright also use  Erd\"{o}s's proof in their number theory text book \cite{HardyWright}.
Let us see what this theorem has to say about the question under investigation.

Let $n$ be an integer with  the diagonal property. That is, given a prime $p$ which does not divide $n$,   $n\, |\, p^{2} - 1$.
Note that if  $n$ divides $p^2 - 1$ then $p^2 - 1 \ge n$, or $p \ge \sqrt{n+1}$.  Equivalently, looking at the contrapositive, we
get the following statement which is more appealing:  If $p < \sqrt{n+1}$, then $p$ divides $n$.

Here is one of several ways to proceed from this point.
Assume that $\sqrt{n+1}/4 \ge 5$ ($\iff n+1 \ge 20^{2}$) and consider the two intervals
\[ \left(\frac{\sqrt{n+1}}{4}, \frac{\sqrt{n+1}}{2}\right), \left(\frac{\sqrt{n+1}}{2}, \sqrt{n+1} \right). \]
 By the Bertrand-Chebyshev theorem each of these intervals has at least one prime. Note that both of these primes are less than $\sqrt{n+1}$. Also, the primes $2, 3$ and $5$ are  less than $\sqrt{n+1}$ because $\sqrt{n+1}/4$ is assumed to be at least $5$. Therefore all these primes, and hence their product, divide $n$. In particular, the product of these primes is at most $n$. From this we have the following inequality
\[ (2)(3)(5) \frac{\sqrt{n+1}}{4} \frac{\sqrt{n+1}}{2} \le n,\]
which simplifies to
\[ 15(n+1) \le 4n.\]
This is impossible.  Therefore we must have $\sqrt{n+1}/4 <  5$, which means $n+1 <  20^{2}$, or $n \le 398$.    Now we claim that  $\sqrt{n+1} \le 7$. If not, then  the product $210$ of the primes $2, 3, 5,$ and $7$ would divide $n$.  Since $n \le 398$, there is only one possibility, namely $n = 210$. But $210$  does not have the
diagonal property because $(11)(191) = 2101 \equiv 1 \mod{210}$. Therefore $\sqrt{n+1} \le 7$ or $n \le 48$. Now let us see what happens if $\sqrt{n+1} > 5$. In this case  the primes $2, 3,$ and $5$ divide $n$. Hence,  their product $30$ divides $n$. The only multiple of $30$ less than $48$ is $30$ itself, which does not have the diagonal property because$(13)(7) = 91 \equiv 1 \mod{30}$. Therefore $\sqrt{n+1} \le 5$, which means $n \le 24$.

The above calculation can be simplified a bit if we use a generalization of the Bertrand-Chebyshev theorem due to Erd\"{o}s, as we will see in the next section.

\section{Theorems of Erd\"{o}s and Ramanujan}
There are several impressive variations and generalizations of the Bertrand-Chebyshev theorem. A generalization due to Ramanujan \cite{Ram1919}, for instance, says that if $n \ge 6$ then there are at least 2 primes between $n$ and $2n$. This theorem was also proved independently later by Erd\"{o}s.  I will use this theorem to simplify the above proof.

Assume that $n$ has the diagonal property. Then, as above, we have the implication, ``$p < \sqrt{n+1} \implies p\, | \, n.$''
Now consider the single interval
\[  \left(\frac{\sqrt{n+1}}{2}, \sqrt{n+1} \right).\]
 If $\sqrt{n+1}/2 \ge 6$ $(\iff$ $n+1\ge 144)$, this interval has at least two primes by  Erd\"{o}s's theorem. Since $\sqrt{n+1}/2 \ge 6$, the primes $2, 3,$ and $5$  will be less than $\sqrt{n+1}$.  Arguing as above, we then have the inequality
\[ (2)(3)(5)\left( \frac{\sqrt{n+1}}{2} \right)^{2} \le  n,\]
which simplifies to give $30(n+1) \le 4n$,  a contradiction. Therefore, $\sqrt{n+1} < 12$, or  $n+1 < 144$.  Now we proceed as before by first showing that $\sqrt{n+1} \le 7$. If not,  then the primes $2, 3, 5$ and $7$, and therefore their product $210$, divides $n$. This is impossible because $n+1 < 144$. Thus  $\sqrt{n+1} \le 7$ which means $n \le 48$.
Similarly if $\sqrt{n+1} > 5$, the primes $2, 3,$ and $5$, and  hence also their product $30$ divides $n$. The only multiple of $30$ less than $48$ is $30$ itself, which does not have the diagonal property. Therefore $\sqrt{n+1} \le 5$, which means $n \le 24$.

\begin{rem}
This is only the beginning.  There are some further generalizations given by Ramanujan \cite{Ram1919}. These  follow right out of his proof of the Bertrand's postulate. To explain these, let $\pi(x)$ denote the number of primes less than or equal to $x$. Ramanujan showed that for each positive integer $k$, there is a prime number $p_k$ such that $\pi(x) - \pi(x/2) \ge k$ if $x \ge p_k$. For example, he showed
\[\pi(x) - \pi(x/2) \ge 1, 2, 3, 4, 5, \cdots \ \ \text{ if } x \ge 2, 11, 17, 29, 41, \cdots \ \ \text{respectively.} \]
The numbers $2, 11, 17, 29, 41, \cdots$ are called the Ramanujan primes. Note that the Bertrand-Chebyshev theorem is covered by the special case
\[\pi(x) - \pi(x/2) \ge 1 \ \ \text{ if } x \ge 2,\]
and the theorem of Erd\"{o}s by the case
\[\pi(x) - \pi(x/2) \ge 2 \ \ \text{ if } x \ge 11.\]
Although one can use these results of Ramanujan to
address our question, the bounds thus obtained become worse and it would take more work to get them down to $24$.
\end{rem}

\section{Further Generalisations}

\begin{center}
\emph{Be Wise! Generalise!}
\end{center}
Instead of working with $\mathbb{Z}_n$, one can look at other rings. For example,  polynomial rings. The question then is: \emph{what are all values of $n$ for which the multiplication table for
$\mathbb{Z}_{n}[x]$ has 1's only on the diagonal?}

Similarly, instead of multiplication tables, one can consider multiplication cubes. This is a natural extension of the notion of a multiplication table and is defined similarly. Given a positive integer $n$, a multiplication cube for $\mathbb{Z}_{n}$ is a cube  $[0, n-1]^{3}$ whose entry at the coordinate $(i, j, k)$  ($0 \le i, j, k \le n-1$) is the product
$ijk \mod{n}$. Now we can ask the same question for these cubes.
\emph{What are all values of $n$ for which the multiplication cube for
$\mathbb{Z}_{n}$ has 1's only on the diagonal?}

I leave it as an amusing project for the reader to venture into these variations and provide  as many different proofs of each  as he or she can.

\vskip 10mm



\end{document}